\documentclass[12pt]{article}
\usepackage[expansion=false]{microtype}
\usepackage[reqno]{amsmath}
\usepackage{fixmath}
\usepackage{amssymb, amsthm, enumitem}
\usepackage{mathtools}
\usepackage{graphicx}
\usepackage{lineno}

\usepackage{ifthen,amssymb,amsfonts,amsmath,amsthm,graphicx, color, tikz, pgf}

\usepackage{hyperref}% wants to be last

\hypersetup{
    plainpages=false,       % needed if Roman numbers in frontpages
    unicode=false,          % non-Latin characters in Acrobat’s bookmarks
    pdftoolbar=true,        % show Acrobat’s toolbar?
    pdfmenubar=true,        % show Acrobat’s menu?
    pdffitwindow=false,     % window fit to page when opened
    pdfstartview={FitH},    % fits the width of the page to the window
    pdfnewwindow=true,      % links in new window
    colorlinks=true,        % false: boxed links; true: colored links
    linkcolor=black,         % color of internal links
    citecolor=black,        % color of links to bibliography
    filecolor=blue,      % color of file links
    urlcolor=blue           % color of external links
}

\newtheoremstyle{plainsl}%
    {\topsep}
    {\topsep}
    {\slshape} % only non-default setting
    {}
    {\normalfont\bfseries}
    {.}
    { }
    {}

\swapnumbers

{\theoremstyle{plainsl}
\newtheorem{thm}{Theorem}[section]
\newtheorem{lem}[thm]{Lemma}

\newtheorem{prop}[thm]{Proposition}
\newtheorem{prob}[thm]{Problem}}
{\theoremstyle{remark}
  \newtheorem{ex}[thm]{Example}
  \newtheorem{rem}[thm]{Remark}
\newtheorem{defn}[thm]{Definition}}

\renewcommand\proof{\noindent\textsl{Proof. }}
\newcommand\sqr[2]{{\vbox{\hrule height.#2pt
   \hbox{\vrule width.#2pt height#1pt \kern#1pt
        \vrule width.#2pt}\hrule height.#2pt}}}
\renewcommand\qed{%
    \ifmmode\eqno\sqr53
    \else\nolinebreak\ \hfill\sqr53\medbreak\fi}

\numberwithin{equation}{section}

\newcommand\pmat[1]{\begin{pmatrix} #1 \end{pmatrix}}

\newcommand{\abs}[1]{\left | #1 \right |}
\newcommand{\lp}{\! \left (}
\newcommand{\rp}{\right )} 
   \newcommand{\lsb}{\left \{ }
\newcommand{\rsb}{\right \} }

\DeclareMathOperator{\Span}{span}
\DeclareMathOperator{\tr}{tr}

\newcommand{\vtxa}{u}
\newcommand{\vtxb}{v}
\newcommand{\vtxc}{w}
\newcommand{\bipartB}{\beta}
\newcommand{\bipartC}{\gamma}
\newcommand{\bipartBC}{\beta \gamma}
\newcommand{\bipartCB}{\gamma \beta}

\definecolor{UWblack}{RGB}{0,0,0}
\definecolor{UWgold}{RGB}{255,213,79}

\title{$P$-Polynomial and Bipartite Coherent Configurations}

\author{Sabrina Lato}

\begin{document}
\maketitle

\begin{abstract}
  We introduce the notion of $P$-polynomial coherent configurations and show that they can have at most two fibres. We then introduce a class of two-fibre coherent configurations which have two distinguished bases for the coherent algebra, similar to the Bose-Mesner algebra of an association scheme. Examples of these bipartite coherent configurations include the $P$-polynomial class of distance-biregular graphs, as well as quasi-symmetric designs and strongly regular designs.
\end{abstract}

\section{Introduction}

Association schemes are a fundamental concept in algebraic combinatorics, linking topics in combinatorial design theory, graph theory, and group theory. In particular, the class of distance-regular graphs arise as \( P \)-polynomial association schemes. The dual notion of \( Q \)-polynomial association schemes was studied by Delsarte~\cite{delsarte1973algebraic} in connection to design theory.

Higman~\cite{higman1970coherent} introduced coherent configurations as a combinatorial generalization of permutation groups. Coherent configurations contain the class of association schemes, which can be seen as combinatorial generalizations of transitive permutation groups. Recently, Suda~\cite{suda2022q} introduced a notion of \( Q \)-polynomial coherent configurations. In this paper, we define \( P \)-polynomial coherent configurations and set up a framework to extend the duality inherent in association schemes to a certain class of coherent configurations.

A motivating example comes from distance-biregular graphs. Distance-biregular graphs are a class of bipartite semiregular graphs that share many combinatorial and algebraic properties with distance-regular graphs. Notable families of distance-biregular graphs include generalized polygons, bipartite distance-regular graphs, partial geometries, and the incidence graphs of certain quasi-symmetric designs.

Godsil and Shawe-Taylor~\cite{distanceRegularised} introduced a notion of locally distance-regular vertices and showed that a graph where every vertex was locally distance-regular is either distance-regular or distance-biregular. In this paper, we take a more algebraic approach and show that the only \( P \)-polynomial coherent configurations are distance-regular and distance-biregular graphs. We then define bipartite coherent configurations, a class of coherent configurations that contains the adjacency algebra of a distance-biregular graph and more closely extends the structure of an association scheme. In addition to distance-biregular graphs, this definition includes quasi-symmetric designs and a class of designs studied by Higman~\cite{higman1988strongly}, Neumaier~\cite{neumaier1980t12}, and Bose, Bridges, Shrikhande, and Singhi~\cite{bose1976characterization, bose1976edge}. We develop some of the properties of bipartite coherent configurations, especially a second distinguished basis for the underlying algebra. We conclude by extending the duality of association schemes to bipartite coherent configurations and highlighting some open problems connecting to the work of Suda~\cite{suda2022q}.

A bipartite coherent configuration can be interpreted as sets of relations between points and blocks, and thus they arise in the contexts of design theory and finite geometry. Restricting to the points or the blocks of a bipartite coherent configuration gives an association scheme, and this perspective has been used in the past to study such structures, for instance, by looking at the strongly regular graphs associated to a quasi-symmetric design. However, bipartite coherent configurations have a deeper structure than just the structure of the underlying association schemes, which is demonstrated in Section~\ref{secKrein}. Thus bipartite coherent configurations give us another way to study incidence structures by extending the theory of association schemes in a natural way.

\section{Coherent Configurations}

Coherent configurations were introduced by Higman~\cite{higman1970coherent}. They are typically defined in terms of sets and relations, though this definition is equivalent to formulation in terms of matrices~\cite{higman1970coherent, higman1975coherent}.

Higman~\cite{higman1987coherent} defined a \textit{coherent algebra} to be a matrix algebra closed under conjugate transpose and Schur multiplication and which contains the all-ones matrix. In other words, a coherent algebra contains both the multiplicative and Schur identities, is closed under both matrix and Schur multiplication, and is also closed under conjugate transpose. Higman~\cite{higman1987coherent} observed that a coherent algebra can be identified by its underlying \textit{coherent configuration}, a family \( \mathcal{C} \) of \( 01 \)-matrices such that:
\begin{enumerate}[label=(\roman*)]
\item Letting \( J \) be the all-ones matrix, \[ \sum_{M \in \mathcal{C}} M = J; \]
\item There exists \( S \subset \mathcal{C} \) such that
  \[ \sum_{M \in S} M = I; \]
\item For all \( M \in \mathcal{C}, \) we have \( M^{T} \in \mathcal{C}; \) and
\item For all \( M_1, M_2 \in \mathcal{C}, \) we have
  \[ M_1 M_2 \in \Span \lp \mathcal{C} \rp. \]
\end{enumerate}

Consequently, we know there exist a set of matrices \( M_1, \ldots, M_f \in \mathcal{C} \) such that
\[ I = M_1 + \cdots + M_f. \]
We may view the rows and columns of the matrices of a coherent configuration as being indexed by some set \( X \). For \( 1 \leq i \leq f, \) we define \( X_i \) to be the subset of \( X \) that is nonzero on the diagonal entries of \( M_i. \) The sets \( X_1, \ldots, X_f \) and the corresponding restriction of relations are the \textit{fibres} of the coherent configuration. Higman~\cite{higman1975coherent} proved that every coherent configuration can be expressed in terms of a standard partition into fibres. We reformulate the matrix theoretic definition of coherent configurations in terms of this standard partition.

Let \( f \geq 1 \) and let \( x_i, t_{ij} \) be positive integers for all \( 1 \leq i, j \leq f. \) Let \( \epsilon_{ij} = 1- \delta_{ij}. \) For \( 1 \leq i, j \leq f \) and \( \epsilon_{ij} \leq r \leq t_{ij} \), let \( A_{r}^{ij} \) be an $x_i \times x_j$ matrix where the entries are either 0 or 1. We define \( \hat{A}_r^{ij} \) to be the \( f \times f \) block matrix where the only nonzero block is \( \lp \hat{A}_r^{ij} \rp_{ij} = A_r^{ij}. \) The set
\[ \mathcal{C} = \lsb \lsb \hat{A}_r^{ij} : \epsilon_{ij} \leq r \leq t_{ij} \rsb : 1 \leq i, j \leq f \rsb \]
is a \textit{coherent configuration} if it satisfies:
\begin{enumerate}[label=(A\arabic*)]
\item \label{cc1}\[ \sum_{i=1}^f A_0^{ii} = I; \]
\item \label{cc2}For all $1 \leq i,j \leq f$, we have
  \[ \sum_{r=~\epsilon_{ij}}^{t_{ij}} A_r^{ij} = J_{x_i,x_j} \]
  where \( J_{x_i,x_j} \) is the \( x_i \times x_j \) matrix of all ones;
\item \label{cc3}For all $1 \leq i, j \leq f$ and $\epsilon_{ij} \leq r \leq t_{ij}$, we have \( \lp \hat{A}_r^{ij} \rp^T \in \mathcal{C}; \) and
\item \label{cc4}For all $1 \leq i, j, h \leq f$ and $\epsilon_{ij} \leq r \leq t_{ij}$ and $\epsilon_{jh} \leq s \leq t_{jh}$, we have
  \[ \hat{A_r}^{ij} \hat{A_s}^{jh} \in \Span \lp \mathcal{C} \rp. \]
\end{enumerate}

\begin{ex}\label{quasisymmetric}Consider a quasi-symmetric 2-\( ( v, k, \lambda ) \) design with block intersection numbers \( s \) and \( t \). This is an incidence structure with \( v \) points and \( b \) blocks where every point is incident to \( r \) blocks, every block is incident to \( k \) points, every pair of points intersect in \( \lambda \) blocks, and every pair of blocks intersect in \( s \) or \( t \) points.

  We can let \( X_0 = I_{v}, \) \( Y_0 = I_b, \) and \( X_1 = J_v-I_v. \) Let \( N_1 \) be the incidence matrix from points to blocks, and let \( N_2 \) be the complement of the incidence matrix. Let \( Y_1 \) be the matrix indexed by the blocks where the \( \lp B, C \rp \)-th entry is 1 if \( \abs{B \cap C} = s \) and 0 otherwise, and let \( Y_2 \) be the matrix with \( \lp B, C \rp \)-th entry equal to 1 if \( \abs{B \cap C} = t. \) The set \\
\resizebox{\textwidth}{\height}{\( \lsb \pmat{X_0 & \mathbf{0} \\ \mathbf{0} & \mathbf{0}}, \pmat{X_1 & \mathbf{0} \\ \mathbf{0} & \mathbf{0}}, \pmat{\mathbf{0} & N_1 \\ \mathbf{0} & \mathbf{0}}, \pmat{\mathbf{0} & N_2 \\ \mathbf{0} & \mathbf{0}}, \pmat{\mathbf{0} & \mathbf{0} \\ N_1^T & \mathbf{0}}, \pmat{\mathbf{0} & \mathbf{0} \\ N_2^T & \mathbf{0}}, \pmat{\mathbf{0} & \mathbf{0} \\ \mathbf{0} & Y_0}, \pmat{\mathbf{0} & \mathbf{0} \\ \mathbf{0} & Y_1}, \pmat{\mathbf{0} & \mathbf{0} \\ \mathbf{0} & Y_2} \rsb \)}
is a coherent configuration.

Conditions~\ref{cc1},~\ref{cc2}, and~\ref{cc3} are obvious. Condition~\ref{cc4} follows from the properties of the incidence matrix of a quasi-symmetric design which can be found, for instance, in Shrikhande and Sane~\cite{shrikhande1991quasi}. Higman~\cite{higman1987coherent} proved that quasi-symmetric designs were coherent configurations and worked out related parameters.
\end{ex}

A coherent configuration \( \mathcal{C} \) forms a basis for the underlying coherent algebra. When considered under Schur multiplication, the elements of \( \mathcal{C} \) are pairwise orthogonal and idempotent. If \( \mathcal{C} \) has a single fibre and the elements of \( \mathcal{C} \) pairwise commute, then \( \mathcal{C} \) forms an association scheme. The coherent algebra of an association scheme is the Bose-Mesner algebra, and has considerably stronger properties than the coherent algebra of an arbitrary coherent configuration~\cite{boseMesner}. Of particular note, the Bose-Mesner algebra has a second distinguished basis where the basis elements are pairwise orthogonal idempotents under ordinary matrix multiplication. Switching between these two bases gives a notion of duality to association schemes.

Given a coherent configuration, Higman~\cite{higman1970coherent, higman1975coherent} defined a second basis of the coherent configuration into orthogonal matrices. However, this choice of second basis is not unique. This led Ito and Munemasa~\cite{ito2020krein} to study the case of fibre-commutative coherent configurations, where they showed the basis was essentially unique. A stronger version of fibre-commutative is fibre-symmetric, satisfying
\begin{enumerate}[resume*]
\item For all $1 \leq i \leq f$ and $0 \leq r \leq t_{ii}$, we have \( \lp A_r^{ii} \rp^T = A_r^{ii}. \)
\end{enumerate}

Suda~\cite{suda2022q} took advantage of the block structure of a coherent configuration to impose even stronger restrictions on the second distinguished basis for a coherent algebra. For \( 1 \leq i, j \leq f, \) let \( \tilde{t}_{ij} = t_{i,j}-\epsilon_{i,j}. \) Suda studied fibre-symmetric coherent configurations where there exist matrices 
\[ \lsb E_{r}^{ij} : 1 \leq i, j \leq f, 0 \leq r \leq \tilde{t}_{ij} \rsb \]
satisfying:
\begin{enumerate}[label=(B\arabic*)]
\item \label{suda1}
\[ E_0^{ij} = \frac{1}{\sqrt{x_i x_j}} J_{x_i, x_j}; \]
\item \label{suda2}The set
\[ \lsb E_r^{ij} : 0 \leq r \leq \tilde{t}_{ij} \rsb \]
is a basis of
\[ \Span \lp \lsb A_r^{ij} : \epsilon_{ij} \leq r \leq t_{ij} \rsb \rp \]
as a vector space;
\item \label{suda3}For all \( 0 \leq r \leq \tilde{t}_{ij}, \) we have $\lp E_r^{ij} \rp^T = E_r^{ji};$ and
\item \label{suda4}For all $1 \leq i, j, j', h \leq f$ and all $0 \leq r \leq \tilde{t}_{ij}$ and $0 \leq s \leq \tilde{t}_{jh}$, we have 
%\begin{equation}\label{198b13}
\[ E_r^{ij} E_s^{j'h} = \delta_{j,j'} \delta_{r,s} E_r^{ih}. \]
%\end{equation}
\end{enumerate}
As we did with \( \hat{A}_r^{ij}, \) it will be convenient to define \( \hat{E}_r^{ij} \) to be the \( f \times f \) block matrix where the only nonzero block is \( \lp \hat{E}_r^{ij} \rp_{ij} = E_r^{ij}. \)

These additional requirements force the coherent configurations to behave spectrally more like an association scheme. We will refer to the set of matrices \( \lsb A_{\epsilon_{ij}}^{ij}, \ldots, A_{t_{ij}}^{ij} \rsb_{i,j=1}^{f} \) as the \textit{Schur idempotents}. With some abuse of notation, we call the matrices \( \lsb E_{0}^{ij}, \ldots, E_{\tilde{t}_{ij}}^{ij} \rsb_{i,j=1}^{f} \) the \textit{spectral idempotents}. These two distinguished bases allow us to develop duality similar to the association scheme case. In particular, Suda~\cite{suda2022q} was studying coherent configurations that were $Q$-polynomial, and we are interested in the dual notion of \( P \)-polynomial coherent configurations.

\section{Distance-Biregular Graphs}\label{secDBG}

A \( P \)-polynomial association scheme is equivalent to the adjacency algebra of a distance-regular graph. Let \( G \) be a graph of diameter \( d \). For \( 0 \leq i \leq d, \) we define \( A_i \) to be the matrix indexed by the vertices of \( G \) such that \( \lp A_i \rp_{\vtxa, \vtxb} = 1 \) if vertices \( \vtxa \) and \( \vtxb \) are at distance \( i \) and 0 otherwise. We will refer to \( A_0, \ldots, A_d \) as the distance adjacency matrices.

The graph \( G \) is distance-regular if for all \( 0 \leq i \leq d, \) there exists a polynomial \( F_i \) of degree \( i \) such that
\[ F_i \lp A_1 \rp = A_i. \]
Distance-regular graphs are a well-studied class of graphs, and more information can be found in the book of Brouwer, Cohen, and Neumaier~\cite{bcn} or the more recent survey of Van Dam, Koolen, and Tanaka~\cite{vanDistance}.

Now suppose that \( G = \lp \bipartB, \bipartC \rp \) is a bipartite graph. Any odd distance matrix will also have the structure of a bipartite graph, so for \( 1 \leq 2i-1 \leq d \) we have a \( \abs{\bipartB} \times \abs{\bipartC} \) matrix \( N_i \) such that
\[ A_{2i-1} = \pmat{\mathbf{0} & N_i \\ N_i^T & \mathbf{0}}. \]
The even distance matrices induce two components, one for each cell in the partition, and so for all \( 0 \leq 2i \leq d \) there exists a \( \abs{\bipartB} \times \abs{\bipartB} \) matrix \( X_i \) and an \( \abs{\bipartC} \times \abs{\bipartC} \) matrix \( Y_i \) such that
\[ A_{2i} = \pmat{X_i & \mathbf{0} \\ \mathbf{0} & Y_i}. \]
We will refer to the \( N_i \) as the \textit{distance-biadjacency matrix} and the \( X_i \) and \( Y_i \) as the distance adjacency matrices for the halved graphs.

\begin{defn}\label{dbg}
The bipartite graph \( G \) is distance-biregular if there exists sequences of polynomials \( \mathcal{P}_0^{\bipartB}, \ldots, \mathcal{P}_{\lfloor \frac{d}{2} \rfloor}^{\bipartB},\mathcal{I}_0^{\bipartB}, \ldots,\mathcal{I}_{\lfloor \frac{d-1}{2} \rfloor}^{\bipartB}, \mathcal{P}_0^{\bipartC}, \ldots,\mathcal{P}_{\lfloor \frac{d}{2} \rfloor}^{\bipartC},\mathcal{I}_0^{\bipartC}, \ldots,\mathcal{I}_{\lfloor \frac{d-1}{2} \rfloor}^{\bipartC} \) such that:
\begin{itemize}
\item For all \( 0 \leq i \leq \lfloor \frac{d}{2} \rfloor, \) the polynomial \( \mathcal{P}_i^{\bipartB} \) has degree \( i \) and satisfies
  \[ \mathcal{P}_i^{\bipartB} \lp N_1 N_1^T \rp = X_i; \]
\item For all \( 1 \leq i \leq \lfloor \frac{d-1}{2} \rfloor, \) the polynomial \(\mathcal{I}_i^{\bipartB} \) has degree \( i \) and satisfies
  \[\mathcal{I}_i^{\bipartB} \lp N_1 N_1^T \rp N_1 = N_i; \]
\item For all \( 0 \leq i \leq \lfloor \frac{d}{2} \rfloor, \) the polynomial \( \mathcal{P}_i^{\bipartC} \) has degree \( i \) and satisfies
  \[ \mathcal{P}_i^{\bipartC} \lp N_1^T N_1 \rp = Y_i; \]
\item For all \( 1 \leq i \leq \lfloor \frac{d-1}{2} \rfloor, \) the polynomial \(\mathcal{I}_i^{\bipartC} \) has degree \( i \) and satisfies
  \[\mathcal{I}_i^{\bipartC} \lp N_1^T N_1 \rp N_1^T = N_i^T. \]
\end{itemize}
\end{defn}

Distance-biregular graphs are closely related to distance-regular graphs. Godsil and Shawe-Taylor~\cite{distanceRegularised} showed they were a natural generalization of distance-regular graphs by proving that any graph where every vertex was locally distance-regular was either distance-regular or distance-biregular. Delorme~\cite{delorme} set up several related algebras, including the adjacency algebra as a coherent configuration. Examples of distance-biregular graphs were given by Delorme~\cite{delorme} and Shawe-Taylor throughout his thesis~\cite{shaweTaylor}, and include generalized polygons, bipartite distance-regular graphs, certain quasi-symmetric designs, partial geometries, and constructions from finite geometry. More recent algebraic formulations have been considered by Fern\'andez and Miklavi\u{c}~\cite{miklavic}, Fern\'andez and Penji\'c~\cite{fernandez2023almost}, and in~\cite{PhD}.

\begin{rem}The definition given here is a reformulation of the definition used in Chapter~3 of~\cite{PhD}. It was shown in Chapter~2 that this definition is equivalent to the ones used by Delorme~\cite{delorme} or Godsil and Shawe-Taylor~\cite{distanceRegularised}.
\end{rem}

Let \( G = \lp \bipartB, \bipartC \rp \) be a distance-biregular graph. Then with the same notation as above, the set
\[\lsb \pmat{X_i & \mathbf{0} \\ \mathbf{0} & \mathbf{0}} \rsb_{i=0}^{\lfloor \frac{d}{2} \rfloor} \bigcup \lsb \pmat{\mathbf{0} & \mathbf{0} \\ \mathbf{0} & Y_j} \rsb_{j=0}^{{\lfloor \frac{d}{2} \rfloor}} \bigcup \lsb \pmat{\mathbf{0} & N_h \\ \mathbf{0} & \mathbf{0}} \rsb_{h=1}^{{\lfloor \frac{d-1}{2} \rfloor}} \bigcup \lsb \pmat{\mathbf{0} & \mathbf{0} \\ N_h^T & \mathbf{0}} \rsb_{h=1}^{{\lfloor \frac{d-1}{2} \rfloor}} \]
forms a coherent configuration. In fact, this coherent configuration is fibre-symmetric and satisfies~\ref{suda1}-~\ref{suda4}. This was established in Chapter~6 of~\cite{PhD}, and we will prove this in a more general setting in Section~\ref{secBasis}.

Distance-regular graphs are \( P \)-polynomial association schemes, and distance-biregular graphs have similar algebraic and combinatorial properties as distance-regular graphs. It is therefore natural to consider distance-biregular graphs as \( P \)-polynomial coherent configurations.

\section{$P$-Polynomial Coherent Configurations}\label{secPPoly}

We are interested in developing a notion of \( P \)-polynomial coherent configurations. We want this definition to include distance-regular and distance-biregular graphs, and also be dual to Suda's definition of \( Q \)-polynomial coherent configurations. In particular, we will assume that we are working with fibre-commutative coherent configurations that satisfy~\ref{suda1} --~\ref{suda4}. From this definition, we will show that every \( P \)-polynomial coherent configuration is equivalent to a distance-regular or distance-biregular graph and thus, in the words of Godsil and Shawe-Taylor~\cite{distanceRegularised}, ``we caught almost exactly the fish we wanted.''

Let \( 1 \leq i, j \leq f. \) Since the Schur idempotents and spectral idempotents form two bases for the same vector space, we know there exist coefficients \( P^{ij}_{r,s} \) such that
\[ A^{ij}_s = \sum_{r=0}^{\tilde{t}_{ij}} P^{ij}_{r,s} E^{ij}_r.\]
The coefficients \( P^{ij}_{r,s} \) are called the \textit{eigenvalues} of the coherent configuration.
Dually, there exist coefficients \( Q^{ij}_{r,s}, \) called the \textit{dual eigenvalues}, satisfying
\[ E^{ij}_s = \sum_{r=\epsilon_{ij}}^{t_{ij}} Q^{ij}_{r,s} A^{ij}_r. \]

Building on the characterization of \( Q \)-polynomial association schemes that can be found, for instance, in Bannai and Ito~\cite{bannai1984algebraic}, Suda~\cite{suda2022q} gave several equivalent definitions of \( Q \)-polynomial association schemes. We will take one of them for our definition.

\begin{defn}\label{qPoly}A coherent configuration is \textit{\( Q \)-polynomial} if, for any \( 1 \leq i, j \leq f, \) there exists a set of polynomials \( \lsb \overline{\nu}_h^{ij} \lp x \rp : 0 \leq h \leq \tilde{t}_{ij} \rsb \) and an ordering of \( E_0^{i,j}, \ldots E_{\tilde{t}_{ij}}^{ij} \) such that for \( 0 \leq h, r \leq \tilde{t}_{ij}, \) the polynomial \( \overline{\nu}_h^{ij} \lp x \rp \) has degree \( h \) and
  \[ Q^{ij}_{h,r} = \overline{\nu}_h^{ij} \lp Q^{ij}_{1,r} \rp. \]
\end{defn}

One way to get a definition of \( P \)-polynomial is rewrite Definition~\ref{qPoly} using the Schur idempotents and the eigenvalues instead of the spectral idempotents and the dual eigenvalues. This definition applied to association schemes immediately gives us the definition of \( P \)-polynomial association schemes that can be found, for instance, in~\cite{bannai1984algebraic}. However, in a distance-biregular graph the eigenvalues of the odd distance matrices are always polynomials of odd degree in the eigenvalues of the adjacency matrix, so the underlying coherent configuration would not satisfy the requirement that the degree of the polynomial corresponds to the ordering of the Schur idempotents.

Thinking of coherent configurations in terms of block matrices, Definition~\ref{qPoly} says there exists a sequence of polynomials for every block. To keep the correlation between the degree of the polynomial and the ordering of the idempotents, we will instead consider a sequence of polynomials one row of the time. 

We will reorder the matrices, so with some abuse of notation, if \( M^{i}_h = A^{ij}_s \) for some $1 \leq j \leq f$ and some $\epsilon_{ij} \leq s \leq t_{ij}$, we will let $P^i_{r,h}$ denote the eigenvalue $P^{ij}_{r,s}$. For \( 1 \leq i \leq f, \) let
\[ t_i := \sum_{j=1}^f t_{ij}. \]

\begin{defn}A coherent configuration is \textit{\( P \)-polynomial} if for all $1 \leq i \leq f$ there is an ordering of the elements of
  \[ \lsb \lsb A^{ij}_h \rsb_{h=\epsilon_{ij}}^{t_{ij}} \rsb_{j=1}^{f} \]
as \( M_0^i, \ldots, M_{t_i}^i \) and a sequence of polynomials \( \nu_0^i, \ldots, \nu_{t_i}^i \) such that, for all \( 0 \leq h \leq t_{i} \) the polynomial $\nu_h^i$ has degree $h$ and satisfies 
\[
\nu_h^i \lp P^i_{r,1} \rp = P^i_{r,h}
\]
 for all \( 0 \leq r \leq \tilde{t}_{i,i}. \)
\end{defn}

Without appealing to distance-biregular graphs, there are linear algebraic ideas to suggest that this is the right definition. Combining~\ref{suda3} and~\ref{suda4}, we see that for \( 1 \leq i, j \leq f \) and \( 0 \leq r \leq \tilde{t}_{ij} \) we have
\[ E_r^{ij} E_r^{ji} = E_r^{ii}. \]
This is the sense in which the spectral idempotents are idempotent, and it suggests a relationship between eigenvalues of different blocks but the same row that is stronger than what exists with the dual eigenvalues.

As a special consequence, we see that \( \tilde{t}_{ii} \geq \tilde{t}_{ij}. \) If \( \tilde{t}_{ij} < \tilde{t}_{ii}, \) it will be convenient to define \( P^{ij}_{r,s} = 0 \) for all \( \tilde{t}_{ij} < r \leq \tilde{t}_{ii}. \) We can also define
\begin{equation}\label{sii826}
  m_r^{i} = \tr \lp E_r^{ii} \rp.
\end{equation}

For the rest of the section, when we have \( P \)-polynomial coherent configuration, we will suppose that \( M_0^{i}, \ldots, M_{t_i}^{i} \) is the ordering of 
  \[ \lsb \lsb A^{ij}_h \rsb_{h=\epsilon_{ij}}^{t_{ij}} \rsb_{j=1}^{f}. \]
Without loss of generality, we may suppose that if \( 0 \leq r < s \leq t_i \) with \( M_r^{i} = A^{ij}_{r'} \) and \( M_s^i = A^{ij}_{s'}, \) then \( r' < s' \).

For a polynomial \( \nu_h^{i} \), the degree tells us which block of the coherent configuration the eigenvalues belong to. This, combined with the structure of matrix multiplication, gives us a lot of information about our coherent configuration.

\begin{lem}\label{oneFibre}If \( \mathcal{C} \) is a \( P \)-polynomial coherent configuration with \( M_1^i = A_1^{ii}, \) then \( \mathcal{C} \) has a single fibre.
\end{lem}

\proof For all \( 0 \leq s \leq \tilde{t}_{ii}, \) let \( \theta_s = P_{s,1}^{ii}. \)

Suppose by induction that there exists some \( 1 \leq r < t_i \) such that \( M_0^{i} = A_0^{ii}, \ldots, M_r^{i} = A_r^{ii}. \) The base case is true by assumption.

By~\ref{cc4}, we know there exist coefficients \( \xi_{1r}^h \) such that
\begin{equation}\label{mg1248}
  A_1^{ii} A_r^{ii} = \sum_{h=0}^{t_{ii}} \xi_{1r}^h A_h = \sum_{h=0}^{t_{ii}} \xi_{1r}^h \sum_{s=0}^{\tilde{t}_{ii}} P^{ii}_{s,h} E_s^{ii}.
\end{equation}
Let \( i_0, \ldots, i_{t_{ii}} \) be the subsequence such that \( M_{i_{h}}^{i} = A_h^{ii}. \) By the inductive hypothesis, \( i_0 = 0, \ldots, i_r = r. \) Then we rewrite Equation~\ref{mg1248} as
\begin{equation}\label{nl1344}
  A_1^{ii} A_r^{ii} = \sum_{h=0}^{t_{ii}} \xi_{1r}^h \sum_{s=0}^{\tilde{t}_{ii}} \nu_{i_h}^{i} \lp \theta_s \rp E_s^{ii}.
\end{equation}

On the other hand, since \( x \nu_r^{i} \lp x \rp \) is a polynomial of degree at most \( t_i, \) we know we can write it as a linear combination of \( \nu_0^{i}, \ldots, \nu_{r+1}^{i}, \) and the coefficient of \( \nu_{r+1}^{i} \) in this expression must be nonzero. Letting \( x= \theta_s, \) we see that \( \theta_s \nu_r^{i} \lp \theta_s \rp \) is a polynomial of degree \( r+1 \) in \( \theta_s. \) Then by~\ref{suda4} and the \( P \)-polynomial property, we have
\begin{equation}\label{nyn2619}
  A_1^{ii} A_r^{ii} = \sum_{s=0}^{\tilde{t}_{ii}} \theta_s P^{ii}_{s,r} E_s^{ii} = \sum_{s=0}^{\tilde{t}_{ii}} \theta_s \nu_r^{ii} \lp \theta_s \rp E_s^{ii}.
\end{equation}
Combining Equation~\ref{nyn2619} with Equation~\ref{nl1344}, we see that \( i_{r+1} = r+1. \)

By induction, we therefore have that \( M_0^{i} = A_0^{ii}, \ldots, M_{t_i}^{i} = A_{t_i}^{ii}, \) from which it immediately follows that \( \mathcal{C} \) can only have one fibre.\qed

Suppose \( \mathcal{C} \) has a single fibre. Since it is clear which fibre we are referring to, we omit the subscripts and superscripts used to refer to the fibre. By the \( P \)-polynomial property, for \( 0 \leq i \leq t, \) we have
\[ A_i = \sum_{r=0}^{t} \nu_i \lp P_{s,1} \rp E_r. \]
By~\ref{suda4} we know that the \( E_r \) are pairwise orthogonal idempotents, so applying the spectral decomposition, which can be found, for instance, in Chapter~8 of Godsil and Royle~\cite{yellow}, we see that
\[ A_i = \nu_i \lp A_1 \rp, \]
which is the definition of a distance-regular graph.

\begin{lem}\label{twoFibre}If \( \mathcal{C} \) is a \( P \)-polynomial coherent configuration with \( M_2^i = A_1^{ii}, \) then \( \mathcal{C} \) has two fibres.
\end{lem}

\proof We know \( M_1^{i} = A_1^{ij} \) for some \( 1 \leq j \leq f \) and \( j \neq i. \) For \( 0 \leq s \leq \tilde{t}_{ij}, \) let \( \theta_s = P^{ij}_{r,1}. \)

Suppose by induction that there exists some \( 1 \leq r \leq \lfloor \frac{t_i}{2} \rfloor -1 \) such that \( M_{2r-1}^{i} = A_{r}^{ij} \) and \( M_{2r}^{i} = A_{r}^{ii}. \) By~\ref{cc3} and~\ref{cc4}, we know that there exist coefficients \( \sigma_{1r}^h \) such that
\[ A_1^{ji} A_r^{ii} = \sum_{h=0}^{t_{ji}} \sigma_{1r}^h A_{h}^{ji} = \sum_{h=0}^{t_{ji}} \sigma_{1r}^h \sum_{s=0}^{\tilde{t}_{ji}} P^{ji}_{s,\ell} E_s^{ii}. \]
Let \( j_1, \ldots, j_{t_{ii}} \) be the subsequence such that \( M_{j_{2r-1}}^{i} = A_r^{ij}. \) Then we rewrite the previous equation as
\begin{equation}\label{yy2243th}
  A_1^{ji} A_r^{ii} = \sum_{h=0}^{t_{ji}} \sigma_{1r}^h \sum_{s=0}^{\tilde{t}_{ji}} \nu_{j_h}^{i} \lp \theta_s \rp E_s^{ji}.
\end{equation}

On the other hand, we can write \( x \nu_{2r}^{i} \lp x \rp \) as a polynomial of degree \( 2r+1. \) By~\ref{suda3},~\ref{suda4}, and the fact that \( \mathcal{C} \) is $P$-polynomial, we have
\begin{equation}\label{ky1304hl}
  A_1^{ji} A_r^{ii} = \sum_{s=0}^{\tilde{t}_{ji}} P^{i}_{s,1} P_{s,2r}^{i} E_s^{ji} = \sum_{s=0}^{\tilde{t}_{ji}} \theta_{s} \nu_{2r}^{i} \lp \theta_s \rp E_s^{ji}.
\end{equation}
Combining Equation~\ref{yy2243th} with~\ref{ky1304hl}, we see that \( j_{r+1} = r+1. \)

A similar argument gives us
\[ A_1^{ij} A_{r+1}^{ji} = \sum_{s=0}^{\tilde{t}_{ji}} \theta_{s} \nu_{2r+1}^{i} \lp \theta_s \rp E_s^{ji} \]
is a polynomial of degree \( 2r+2 \) in \( \theta_s, \) and
\[ A_1^{ij} A_{r+1}^{ji} = \sum_{h=0}^{t_{ii}} \xi_{1 \lp r+1 \rp}^h \sum_{s=0}^{\tilde{t}_{ii}} \nu_{i_h}^{i} \lp \theta_s \rp E_s^{ii}, \]
so \( M_{2r+2}^{i} = A_{r+1}^{ii}. \)

By induction, for all \( 0 \leq h \leq t_i, \) we must have \( M_h^{i} = A_{h/2}^{ii} \) if \( h \) is even and \( A_{(h-1)/2}^{ij} \) if \( h \) is odd. It follows that \( \mathcal{C} \) contains two fibres.\qed

Now suppose that \( \mathcal{C} \) is \( P \)-polynomial with two fibres. It must have the form just described, where the polynomials \( \nu_{2h}^{1} \) and \( \nu_{2h}^{2} \) are even and correspond to the eigenvalues of \( A^{11} \) and \( A^{22}, \) and the polynomials \( \nu_{2h+1}^{1} \) are odd and correspond to eigenvalues of \( A^{12}. \)

Let \( 0 \leq 2h \leq t_1. \) Since \( \nu_{2h}^{1} \) is even,~\ref{suda1} --~\ref{suda4} tell us that
\[ A^{11}_h = \sum_{r=0}^{\tilde{t}_{11}} \nu_{2h}^{1} \lp P^{12}_{r,1} \rp E^{11}_r = \sum_{r=0}^{\tilde{t}_{12}} \nu_{2h}^{1} \lp P_{r,1}^{12} \rp E_r^{12} E_r^{21} + \nu_{2h}^{1} \lp P_{\tilde{t}_{11}}^{12} \rp \lp I - \sum_{r=0}^{\tilde{t}_{12}} E_r^{12} E_r^{21} \rp \]
is a polynomial in \( NN^T. \) Thus \( \nu_0^{1}, \nu_2^{1}, \ldots, \nu_{2t_{11}}^{1} \) gives us the sequence \( \mathcal{P}_0^{\bipartB}, \ldots, \mathcal{P}_{t_{11}}^{\bipartB} \) of Definition~\ref{dbg}. Similarly, for all \( 0 \leq 2h \leq t_2, \) the matrix \( A^{22}_h \) is a polynomial in \( N^TN \) giving us the sequence \( \mathcal{P}_0^{\bipartC}, \ldots, \mathcal{P}_{t_{22}}^{\bipartC}. \)

We can write any odd polynomial \( \nu_{2h+1} \lp x \rp \) as an even polynomial multiplied by \( x \). Then we see that for \( 1 \leq 2h+1 \leq t_1, \) there exists a polynomial \( \mathcal{I}_h^{\bipartB} \) such that we can write
\[ A_h^{12} = \sum_{r=0}^{\tilde{t}_{11}} \nu_{2h+1}^{1} \lp P^{12}_{r,1} \rp E_r^{12} = \mathcal{I}_h^{\bipartB} \lp NN^T \rp N. \]
An analogous argument holds for \( A_h^{21}, \) so by Definition~\ref{dbg} we see \( \mathcal{C} \) is the adjacency algebra of a distance-regular graph.

This leads us our main result.

\begin{thm}\label{redFibre}If \( \mathcal{C} \) is a \( P \)-polynomial coherent configuration, then it is the adjacency algebra of a distance-regular or distance-biregular graph.
\end{thm}

\proof Let \( M_0^{i}, \ldots, M_{t_i}^i \) be the ordering, and let \( 1 \leq j \leq f \) such that that \( M_1^{i} = A_1^{ij}. \)

By~\ref{suda3} and~\ref{suda4}, we know
\begin{equation}\label{yi5438ic}
 A_1^{ij} A_1^{ji} = \sum_{\ell=0}^{\tilde{t}_ii} \lp P^{ij}_{r,1} \rp^2 E_r^{ii}.
\end{equation}
Let \( i_0, \ldots, i_{t_{ii}} \) be the subsequence of \( 0, \ldots, t_i \) such that \( M_{i_{h}}^{i} = A_{h}^{ii}. \) Then since \( \mathcal{C} \) is $P$-polynomial, we know there exist coefficients \( \xi_{11}^h \) such that
\begin{equation}\label{bt5454ww}
A_1^{ij} A_1^{ji} = \sum_{h=0}^{t_{ii}} \xi_{11}^h \sum_{r=0}^{\tilde{t}_ii} \nu_{i_{h}}^i \lp P^{ii}_{r,1} \rp E_r^{ii}.
\end{equation}

Combining Equations~\ref{yi5438ic} and~\ref{bt5454ww}, we see that \( \xi_{11}^h = 0 \) unless \( h \leq 2. \) Further, we must have that \( 2 \in \lsb i_1, i_2 \rsb. \) If \( i_1 =2, \) then Lemma~\ref{twoFibre} tells us that \( \mathcal{C} \) is the adjacency algebra of a distance-biregular graph, and if \( i_2 =2, \) then Lemma~\ref{oneFibre} tells us that \( \mathcal{C} \) is the adjacency algebra of a distance-regular graph.\qed

\section{Bipartite Coherent Configurations}\label{secBCC}

For the rest of the paper, we are interested in the class of coherent configurations that contain distance-biregular graphs. By restricting ourselves to two fibre coherent configurations that are fibre symmetric and satisfy~\ref{suda1}-~\ref{suda4}, we can set up notions of duality that parallel what is happening with association schemes. In fact, we will define a particular class of coherent configurations in terms of the Schur idempotents, and use this definition to develop a spectral basis satisfying~\ref{suda1}-~\ref{suda4}. Since we are now dealing with a small and fixed number of fibres, we can also reformulate the definitions to avoid some of the subscripts and superscripts.

Let \( \bipartB \) and \( \bipartC \) be sets and let \(  t_{\bipartB}, t_{\bipartC}, \) and \( t_{\bipartBC} \) be positive integers. Let \( \lsb X_0, \ldots, X_{t_{\bipartB}} \rsb \) be a set of \( \abs{\bipartB} \times \abs{\bipartB} \) 01-matrices. Similarly, let \( \lsb Y_0, \ldots, Y_{t_{\bipartC}} \rsb \) be a set of \( \abs{\bipartC} \times \abs{\bipartC} \) 01-matrices, and let \( \lsb N_1, \ldots, N_{t_{\bipartBC}} \rsb \) be a set of \( \abs{\bipartB} \times \abs{\bipartC} \) 01-matrices. Let
\[ \mathcal{C} = \lsb \pmat{X_i & \mathbf{0} \\ \mathbf{0} & \mathbf{0}} \rsb_{i=0}^{t_{\bipartB}} \bigcup \lsb \pmat{\mathbf{0} & \mathbf{0} \\ \mathbf{0} & Y_j} \rsb_{j=0}^{t_{\bipartC}} \bigcup \lsb \pmat{\mathbf{0} & N_h \\ \mathbf{0} & \mathbf{0}} \rsb_{h=1}^{t_{\bipartBC}} \bigcup \lsb \pmat{\mathbf{0} & \mathbf{0} \\ N_h^T & \mathbf{0}} \rsb_{h=1}^{t_{\bipartBC}}. \]
We say that $\mathcal{C}$ forms a \textit{bipartite coherent configuration} if it satisfies the following:
\begin{enumerate}[label=(C\arabic*)]
\item\label{bcc1} The matrices \( X_0 \) and \( Y_0 \) are the identity matrix;
\item\label{bcc2}
  \[ \sum_{M \in \mathcal{C}} M = J; \]
\item\label{bcc3} For all \( M \in \mathcal{C}, \) we have \( M^T \in \mathcal{C}. \)
\item\label{bcc5} For all \( M_1, M_2 \in \mathcal{C}, \) we have \( M_1M_2 \in \Span \lp \mathcal{C} \rp; \)
\item\label{bcc4} For all $1 \leq i,j \leq t_{\bipartBC}$, we have
\[ N_iN_j^T = N_jN_i^T \]
and
\[ N_i^T N_j = N_j^T N_i; \]
\item\label{bcc6} The set
\[ \lsb N_i N_j^T : 1 \leq i, j \leq t_{\bipartBC} \rsb \cup \lsb I \rsb  \]
spans \( \Span \lsb X_i : 0 \leq i \leq t_{\bipartB} \rsb \) and
\[ \lsb N_i^TN_j : 1 \leq i, j \leq t_{\bipartBC} \rsb \cup \lsb I \rsb \]
spans \( \Span \lsb Y_i : 0 \leq i \leq t_{\bipartC} \rsb. \)
\end{enumerate}

\begin{rem}It may be possible to weaken assumption~\ref{bcc4} to \( N_iN_j^T \) and \( N_i^T N_j \) are normal to give a closer analogue to association schemes that are commutative, but not necessarily symmetric, then the definition given here.
\end{rem}

The adjacency algebra of a distance-biregular graph is a bipartite coherent configuration. That perspective was made explicit in Chapter~6 of~\cite{PhD}, but the structure of bipartite coherent configurations was used implicitly throughout the rest of the thesis.

Higman~\cite{higman1987coherent} defined the \textit{type} of a coherent configuration, denoted by a matrix indexed by the fibres where the entries indicate the number of elements in the corresponding block of the coherent configuration. For the purposes of two-fibre coherent configurations, we can think of them as symmetric matrices, denoted
\[ \pmat{t_{\bipartB} +1 & t_{\bipartBC} \\ & t_{\bipartC}+1}. \]
%or \( \lp t_{\bipartB}+1, t_{\bipartBC}; t_{\bipartC}+1 \rp. \)

Quasi-symmetric designs are coherent configurations of type \( \pmat{2 & 2 \\ & 3}, \)
and in fact, they are the only coherent configurations of that type~\cite{higman1987coherent}. We can consider some other coherent configurations of small type in the context of bipartite coherent configurations.

\begin{ex}\label{stronglyRegularDesign}Higman~\cite{higman1988strongly} termed the coherent configurations of type \( \pmat{3 & 2 \\ & 3} \) \textit{strongly regular designs}. As before, we let \( v \) be the number of points and \( b \) be the number of blocks, let \( k \) be the number of points incident to any given block and let \( r \) be the number of blocks incident to a given point. Strongly regular designs can be phrased combinatorially in a number of ways, but expressed in terms of matrices of a coherent configuration we have numbers \( s_1 >  s_2 \) and \( t_1 > t_2 \) such that
  \[ N_1 N_1^T = rI + s_1 X_1 + s_2 X_2. \]
and
  \begin{equation}\label{srd}
    N_1^T N_1 = kI + t_1 Y_1 + t_2 Y_2.
  \end{equation}
  
  By~\ref{cc2}, we have that
  \begin{equation}\label{srd2}
    N_1^T N_2 = N_1^T \lp J_{v,b}- N_1 \rp = kJ_{b,b} - N_1^T N_1 = \lp J_{v,b}-N_1 \rp^T N_1 = N_2^T N_1
  \end{equation}
  establishing~\ref{bcc4}. Combining Equations~\ref{srd} and~\ref{srd2} we have
  \[ N_1^T N_2 = \lp k- t_1 \rp Y_1 + \lp k-t_2 \rp Y_2, \]
  so
  \begin{equation}\label{srd3}
    Y_1 + Y_2  = \frac{1}{k} \lp N_1^TN_2 + N_1^T N_1 -kI \rp.
  \end{equation}
  Now combining Equations~\ref{srd} and~\ref{srd3} we have
  \begin{equation}\label{srd4}
    \lp t_1 - t_2 \rp Y_2 = \frac{t_1}{k} \lp N_1^TN_1 + N_1^T N_1 - kI \rp -N_1^T N_1 + kI
  \end{equation}
  and
  \begin{equation}\label{srd5}
    \lp t_2 - t_1 \rp Y_1 = \frac{t_2}{k} \lp N_1^TN_1 + N_1^T N_1 - kI \rp -N_1^T N_1 + kI,
  \end{equation}
  which gives us~\ref{bcc6}. The same argument holds for \( X_1 \) and \( X_2, \) which shows us that a strongly regular design is indeed a bipartite coherent configuration.

Higman~\cite{higman1988strongly} worked out further parameters of strongly regular designs, and Hanaki~\cite{hanaki2021note} expanded on and corrected some of the proofs. Strongly regular designs have also been studied under other names.

  Neumaier~\cite{neumaier1980t12} studied the notion of a \( t \frac{1}{2} \)-design, with a structure between that of a \( t \)-design and a \( \lp t+1 \rp \)-design. The \( 1 \frac{1}{2} \)-designs contain strongly regular designs. Neumaier further proved that every \( 2 \frac{1}{2} \)-design is either a 3-design or a multiple of a symmetric 2-design, every \( 3 \frac{1}{2} \)-design is either a 4-design or a multiple of a Hadamard 3-design, and every \( t \frac{1}{2} \)-design for \( t \geq 4 \) is a \( \lp t+1 \rp \)-design. Therefore, the strongly regular designs are the most interesting case of \( t \frac{1}{2} \)-designs.

Bose, Shrikhande, and Singhi~\cite{bose1976edge} studied a multigraph extension of strongly regular graphs, which they called partial geometric designs, as a multigraph analogue to partial geometries. Bose, Bridges, and Shrikhande~\cite{bose1976characterization} further explored the spectral properties of partial geometric designs, which are equivalent to the \( 1 \frac{1}{2} \) designs of Neumaier~\cite{neumaier1980t12} and contain strongly regular designs.
\end{ex}

It is also worth noting that Equations~\ref{srd4} and~\ref{srd5} are precisely the equations given by a quasi-symmetric design, so quasi-symmetric designs give us another large family of bipartite coherent configurations. However, not every coherent configuration with two fibres is a bipartite coherent configuration.

\begin{ex}\label{hobart}Hobart~\cite{hobartPhd} studied coherent configurations of type \[ \pmat{2 & 2 \\ & 4}. \] These are the coherent configurations where there exist coefficients \( t_1, t_2, t_3 \) such that
%    \begin{equation}\label{ttf1}
    \[ N_1^T N_1 = kI + t_1 Y_1 + t_2 Y_2 + t_3 Y_3. \]
%    \end{equation}

    By~\ref{cc2}, we have that
    \[ N_1^T N_2 = N_1^T \lp J-N_1 \rp = kJ - N_1^T N_1 = \lp J-N_1 \rp^T N_1, \]
    so this coherent configuration does satisfy~\ref{bcc4}. However, we have
    
    %    \begin{equation}\label{ttf2}
    \[ N_1^T N_2 = \lp k-t_1 \rp Y_1 + \lp k-t_2 \rp Y_2 + \lp k - t_3 \rp Y_3, \]
%    \end{equation}
    and
    %    \begin{equation}\label{ttf3}
    \[ N_2^T N_2 = \lp b-k \rp I + \lp b-2k+t_1 \rp Y_1 + \lp b- 2k + t_2 \rp Y_2 + \lp b-2k+t_3 \rp Y_3. \]
%    \end{equation}

    Combining these, we get
    \[ N_2^T N_2 = \frac{\lp b-2k \rp}{k} \lp N_1^TN_1 + N_1^T N_2 \rp + N_1^TN_1, \]
    so the span of \( \lsb N_1N_1^T, N_1^TN_2, N_2^T N_2 \rsb \) has dimension two, which contradicts~\ref{bcc6}. Thus the two-fibre coherent configurations considered by Hobart~\cite{hobartPhd} are not bipartite coherent configurations.  
\end{ex}

\section{A Second Basis}\label{secBasis}

The following is a standard fact in linear algebra that can be found, for instance, in Section~10.~3 of Godsil~\cite{blue}.

\begin{lem}[~\cite{blue}]\label{BBT}Let \( N \) be a matrix. Then the nonzero eigenvalues of \( NN^T \) and \( N^TN \) are the same with the same multiplicity.
\end{lem}

Consider a set $S$ of $n \times n$ commuting idempotent matrices. Let $F, E \in S$.  If $FE = E,$ we say that $E \leq F$. This relation defines a partial order on $S$. We say that $E \in S$ is \emph{minimal} if there does not exist another matrix $F \in S$ such that $F \leq E$. The following result, which can be found in Section~3.~4 of Godsil and Meagher~\cite{grey}, will be useful.

\begin{thm}[Godsil and Meagher~\cite{grey}]\label{basisIdempotent}
  Let  $\mathcal{S} = \lsb A_0, \ldots, A_d \rsb$ be a set of normal, commuting matrices. There exists a set of minimal idempotents $\lsb E_0, \ldots, E_t \rsb$ that form a basis for \( \Span \lp \mathcal{S} \rp \) and are pairwise orthogonal. Further, these minimal idempotents are projections into the eigenspaces of \( \mathcal{S}. \)
\end{thm}

We can use this set of minimal idempotents to construct a second basis for our bipartite coherent configuration.

\begin{thm}\label{dualBasis} Let \( \mathcal{C} \) be a bipartite coherent configuration. Then \( \Span \lp \mathcal{C} \rp \) has a second basis
  \[ \mathcal{B} = \lsb \pmat{L_i & \mathbf{0} \\ \mathbf{0} & \mathbf{0}} \rsb_{i=0}^{t_{\bipartB}} \bigcup \lsb \pmat{\mathbf{0} & \mathbf{0} \\ \mathbf{0} & R_j} \rsb_{j=0}^{t_{\bipartC}} \bigcup \lsb \pmat{\mathbf{0} & D_h \\ \mathbf{0} & \mathbf{0}} \rsb_{h=0}^{\tilde{t}_{\bipartBC}} \bigcup \lsb \pmat{\mathbf{0} & \mathbf{0} \\ D_h^T & \mathbf{0}} \rsb_{h=0}^{\tilde{t}_{\bipartBC}} \]  
satisfying
\begin{enumerate}[label=(D\arabic*)]
\item\label{dbc3} For all \( 0 \leq r,s \leq t_{\bipartB}, \) we have
  \[ L_r L_s = \delta_{rs} L_r \]
  and similarly for all \( 0 \leq r, s \leq t_{\bipartC}, \) we have
  \[ R_r R_s = \delta_{rs} R_r; \]
\item\label{dbc1} \( L_0 = \frac{1}{\abs{\bipartB}} J_{\bipartB} \), \( D_0 = \frac{1}{\sqrt{\abs{\bipartB} \abs{\bipartC}}} J_{\bipartB, \bipartC}, \) and \( R_0 = \frac{1}{\abs{\bipartC}} J_{\bipartC}; \)
\item\label{dbc2} \[ \sum_{r=0}^{t_{\bipartB}} L_r = I \] and \[ \sum_{r=0}^{t_{\bipartC}} R_r = I; \]
\item\label{dbc4} For all \( 0 \leq r \leq \tilde{t}_{\bipartBC}, \) we have \( L_r = D_r D_r^T \) and \( R_r = D_r^T D_r; \) and
\item\label{dbc5} For all $M_1, M_2 \in \mathcal{B}$, we have $M_1 \circ M_2 \in \mathcal{B}$.
\end{enumerate}
\end{thm}

\proof Let 
\[ \mathcal{S} = \bigcup_{i=1}^{t_{\bipartB \bipartC}} \lsb \pmat{\mathbf{0} & N_i \\ N_i^T & \mathbf{0}} \rsb \bigcup \lsb \pmat{I_{\bipartB} & \mathbf{0} \\ \mathbf{0} & I_{\bipartC}} \rsb. \]
Clearly \( \Span \lp \mathcal{S} \rp  \subseteq \Span \lp \mathcal{C} \rp \), and by~\ref{bcc4} we know that the elements of \( \mathcal{S} \) pairwise commute. Then by Theorem~\ref{basisIdempotent}, we know that \( \Span \lp \mathcal{S} \rp \) has a basis of minimal idempotents which represent orthogonal projection into the eigenspaces. We call this basis $\mathcal{T}$.

Let
\[ F = \frac{1}{2} \pmat{\frac{1}{\abs{\bipartB}} J_{\bipartB} & \frac{1}{\sqrt{\abs{\bipartB}\abs{\bipartC}}} J_{\bipartB, \bipartC} \\ \frac{1}{\sqrt{\abs{\bipartB}\abs{\bipartC}}} J_{\bipartC, \bipartB} & \frac{1}{\abs{\bipartC}} J_{\bipartC}}. \]
Note that \( F \) is idempotent, and by~\ref{bcc2} we have
\[ \sum_{i=1}^{t_{\bipartBC}} \pmat{\mathbf{0} & N_i \\ N_i^T & \mathbf{0}} F = \pmat{\mathbf{0} & J_{\bipartB, \bipartC} \\ J_{\bipartC, \bipartB} & \mathbf{0}} F = \sqrt{\abs{\bipartB} \abs{\bipartC}} F, \]
Since \( F \) has rank one, it must be a minimal idempotent, so \( F \in \mathcal{T}. \)

Let \( E_r \in \mathcal{T}. \) We can break $E_r$ into the same sized blocks that we used for the matrices in \( \mathcal{C} \). That is, there exist matrices \( L_r, D_r, \) and \( R_r \) where \( L_r \) has rows and columns indexed by \( \bipartB, \) \( D_r \) has rows indexed by \( \bipartB \) and columns by \( \bipartC, \) and \( R_r \) has rows and columns indexed by \( \bipartC, \) such that
\[ E_r = \frac{1}{2} \pmat{L_r & \mathbf{0} \\ \mathbf{0} & \mathbf{0}} + \frac{1}{2} \pmat{\mathbf{0} & D_r \\ \mathbf{0} & \mathbf{0}} + \frac{1}{2} \pmat{\mathbf{0} & \mathbf{0} \\ D_r^T & \mathbf{0}} + \frac{1}{2} \pmat{\mathbf{0} & \mathbf{0} \\ \mathbf{0} & R_r}. \]

Condition~\ref{dbc1} follows immediately from the fact that \( F \in \mathcal{T}. \)

Suppose that \( E_r \in \mathcal{T}, \) and let
\[ \tilde{E}_r = \frac{1}{2} \pmat{L_r & \mathbf{0} \\ \mathbf{0} & \mathbf{0}} + \frac{1}{2} \pmat{\mathbf{0} & -D_r \\ \mathbf{0} & \mathbf{0}} + \frac{1}{2} \pmat{\mathbf{0} & \mathbf{0} \\ -D_r^T & \mathbf{0}} + \frac{1}{2} \pmat{\mathbf{0} & \mathbf{0} \\ \mathbf{0} & R_r}. \]
We claim that \( \tilde{E}_r \in \mathcal{T}. \)

Recall that \( E_r \) represents projection into the \( \theta_r \) eigenspace of some \( A_i \in S. \) Since \( A_i \) can be seen as the adjacency matrix of a bipartite graph, a simple argument that can be found, for instance, in Section~8.~8 of Godsil and Royle~\cite{yellow}, tells us that \( \tilde{E}_r \) represents projection into the \( -\theta_r \) eigenspace of \( A_i. \) Thus \( \tilde{E}_r \in \Span \lp \mathcal{T} \rp, \) so by the minimality of elements of \( \mathcal{T} \) there exists some nonzero matrix \( E_s \) satisfying
\[ E_s = \frac{1}{2} \pmat{L_s & D_s \\ D_s^T & R_s} = E_s \tilde{E}_r = \frac{1}{4} \pmat{L_sL_r - D_sD_r^T & -L_sD_r + D_s R_r \\ D_s^T L_s -R_s D_r^T & -D_s^T D_r + R_sR_r}. \]
But then \( \tilde{E}_s \in \Span \lp \mathcal{T} \rp \) and satisfies 
\[  \tilde{E}_s E_r = \pmat{L_s L_r - D_sD_r^T & L_s D_r-D_s R_r \\ -D_s^T L_r + R_s D_r & -D_s^T D_r+R_s R_r} = \pmat{L_s & -D_s \\ -D_s^T & L_s}, \]
and thus $\tilde{E}_s \leq E_r$. Given the choice of $E_r$ as minimal, this means that $\tilde{E}_s = E_r$, and thus  \( \tilde{E}_r = E_s \in \mathcal{T}. \)

This allows us to split \( \mathcal{T} \) into a set \( \mathcal{T}' \subseteq \mathcal{T} \) such that
\[ \bigcup_{E \in \mathcal{T'}} \lsb E, \tilde{E} \rsb = \mathcal{T} \]
Let
\[ \mathcal{B} = \bigcup_{E_r \in \mathcal{T}'} \lsb \pmat{L_r & \mathbf{0} \\ \mathbf{0} & \mathbf{0}}, \pmat{\mathbf{0} & D_r \\ \mathbf{0} & \mathbf{0}}, \pmat{\mathbf{0} & \mathbf{0} \\ D_r^T & \mathbf{0}}, \pmat{\mathbf{0} & \mathbf{0} \\ \mathbf{0} & R_r} \rsb. \]
Clearly the elements of \( \mathcal{B} \) span \( \mathcal{S}. \) We wish to show that, removing any zero matrices, they form a basis for \( \Span \lp \mathcal{C} \rp. \) With some abuse of notation, we refer to the matrices \( L_r, D_r, \) and \( R_r \) as being in \( \mathcal{B}. \)

Note that
\[ E_r + \tilde{E}_r = \pmat{L_r & \mathbf{0} \\ \mathbf{0} & R_r} \]
is idempotent, from which we see that $L_r$ and $R_r$ are themselves idempotent. If \( E_r, E_s \in \mathcal{T} \) with $s \neq r$, then
\[ \mathbf{0} = \lp E_r + \tilde{E}_r \rp \lp E_s + \tilde{E}_s \rp= \pmat{L_r & \mathbf{0} \\ \mathbf{0} & R_r} \pmat{L_s & \mathbf{0} \\ \mathbf{0} & R_s}, \]
so the $L_r$ are pairwise orthogonal, as are the $R_r$. This establishes~\ref{dbc3}.

In a similar vein, if \( D_r, D_s \) are nonzero matrices in \( \mathcal{B}, \) then 
\[ \pmat{\mathbf{0} & D_r \\ D_r^T & \mathbf{0}} \pmat{\mathbf{0} & D_s \\ D_s^T & \mathbf{0}} = \lp E_r - \tilde{E}_r \rp \lp E_s - \tilde{E}_s \rp = \delta_{rs} \pmat{L_r & \mathbf{0} \\ \mathbf{0} & R_r}. \]
Reordering the matrices as necessary, this establishes~\ref{dbc4}.

We have
\[ \sum_{L_r \in \mathcal{B}} \pmat{L_r & \mathbf{0} \\ \mathbf{0} & \mathbf{0}} + \sum_{R_r \in \mathcal{B}} \pmat{\mathbf{0} & \mathbf{0} \\ \mathbf{0} & R_r} = \sum_{E_r \in \mathcal{T}} E_r = \pmat{I_{\bipartB} & \mathbf{0} \\ \mathbf{0} & I_{\bipartC}}, \]
establishing~\ref{dbc2}.

To show \( \mathcal{B} \) is a basis for  \( \Span \lp \mathcal{C} \rp, \) we use~\ref{bcc6}. Let \( 0 \leq i, j \leq \tilde{t}_{\bipartBC}. \) Since \( \mathcal{B} \) spans \( \mathcal{S}, \) we know there exist scalars \( P_{r,i}^{\bipartBC} \) such that
\[ N_i = \sum_{D_r \in \mathcal{B}} P_{r,i}^{\bipartBC} D_r, \]
and similarly there exist scalars \( P_{s,j}^{\bipartBC} \) such that
\[ N_j^T = \sum_{D_s^T \in \mathcal{B}} P_{s,j}^{\bipartBC} D_s^T, \]
so
\[ N_iN_j^T = \sum_{D_r \in \mathcal{B}} \sum_{D_s^T \in \mathcal{B}} P_{r,i}^{\bipartBC} P_{s,j}^{\bipartBC} D_r D_s^T = \sum_{D_r \in \mathcal{B}} P_{r,i}^{\bipartBC} P_{s,i}^{\bipartBC} L_r. \]
An analogous argument holds for \( N_i^T N_j \), so by~\ref{bcc6} and~\ref{dbc4}, we see that \( \mathcal{B} \) spans \( \mathcal{C} \).

Finally, the matrices in \( \mathcal{C} \) are idempotent under Schur multiplication, so~\ref{dbc5} follows from the fact that \( \mathcal{B} \) and \( \mathcal{C} \) are both bases for the same space.\qed

Note that~\ref{suda2} and~\ref{suda3} are implicit in the description of \( \mathcal{B}. \) Further,~\ref{suda1} is equivalent to~\ref{dbc1}, and~\ref{dbc3} and~\ref{dbc4} combine to be equivalent to~\ref{suda4}, so bipartite coherent configurations are indeed examples of the types of coherent configurations Suda~\cite{suda2022q} considered.

\section{Eigenvalues and Dual Eigenvalues}

In this section, we lay out some of the properties of the eigenvalues and dual eigenvalues for bipartite coherent configurations. This has been done more generally~\cite{higman1970coherent,higman1987coherent,suda2022q} and more specifically~\cite{bose1976characterization,bcn,PhD} by a number of authors in various contexts, but it is useful to fix some notation and work within the framework of bipartite coherent configurations.

Recall that
\begin{equation}
N_i = \sum_{j=0}^{\tilde{t}_{\bipartBC}} P_{j,i}^{\bipartBC} D_j. \label{eq:1aecc5}
\end{equation}
gives us the eigenvalues associated to \( \bipartBC. \) We define the \textit{eigenmatrix} for $\bipartBC$ as the $t_{\bipartBC} \times t_{\bipartBC}$ matrix $P^{\bipartBC}$ with \(\lp j, i-1 \rp \)-th entry equal to \( P^{\bipartBC}_{j,i}. \) The eigenvalues associated to \( \bipartBC \) are also the eigenvalues associated to \( \bipartCB. \) 

For \( 1 \leq i \leq t_{\bipartBC}, \) the \textit{valency associated to $\bipartBC$}, denoted $k_i^{\bipartBC}$ is the row sum of $N_i.$ From~\ref{dbc1} we see that
\[ N_j D_0^T = P_{j,0}^{\bipartBC} L_0 = \frac{1}{\abs{\bipartB}} P_{j,0}^{\bipartBC} J_{\bipartB}, \]
so \( N_i \) has constant row sum
\begin{equation}\label{eq:kl2014}
  \sqrt{\abs{\bipartC}/\abs{\bipartB}} P^{\bipartBC}_{0,i}.
\end{equation}
We can similarly define the valency associated to \( \bipartCB \) as the row sum of \( N_i^T, \) and the same argument gives us
\[ k_i^{\bipartCB} = \sqrt{\abs{\bipartB}/\abs{\bipartC}} P^{\bipartBC}_{0,i}.\]

Dually to Equation~\ref{eq:1aecc5}, we have
\begin{equation}
D_i = \frac{1}{\sqrt{\abs{\bipartB} \abs{\bipartC}}} \sum_{j=1}^{t_{\bipartBC}} Q_{j,i}^{\bipartBC} N_j \label{eq:68c2ad}
\end{equation}
and the \textit{dual eigenmatrix} for \( \bipartBC \) is the \( t_{\bipartBC} \times t_{\bipartBC} \) matrix with \( \lp j-1, i \rp \)-th entry equal to \( Q_{j,i}^{\bipartBC} \).

By~\ref{dbc1} and~\ref{bcc2}, we have
\[ D_0 = \frac{1}{\sqrt{\abs{\bipartB} \abs{\bipartC}}} J_{\bipartB, \bipartC} = \frac{1}{\sqrt{\abs{\bipartB} \abs{\bipartC}}} \sum_{j=1}^{t_{\bipartBC}} Q_{j,0}^{\bipartBC} N_j, \]
and thus \( Q_{i0}^{\bipartBC} = 1 \) for all \( 1 \leq j \leq t_{\bipartBC}. \)

Let \( 1 \leq i \leq t_{\bipartBC}. \) We have
\[ N_i = \sum_{h=0}^{\tilde{t}_{\bipartBC}} P^{\bipartBC}_{ih} D_h  = \frac{1}{\sqrt{\abs{\bipartB} \abs{\bipartC}}} \sum_{h=0}^{\tilde{t}_{\bipartBC}} \sum_{j=1}^{t_{\bipartBC}} P_{ih}^{\bipartBC} Q_{hj}^{\bipartBC} N_j, \]
from which it follows that, for all \( 1 \leq i, j \leq t_{\bipartBC} \) we must have
\[ \sqrt{\abs{\bipartB} \abs{\bipartC}} \delta_{ij} = \sum_{h=0}^{\tilde{t}_{\bipartBC}} P^{\bipartBC}_{ih} Q^{\bipartBC}_{hj} = \lp P^{\bipartBC} Q^{\bipartBC} \rp_{ij}, \]
so
\[ P^{\bipartBC} Q^{\bipartBC} = \sqrt{\abs{\bipartB}\abs{\bipartC}}I_{t_{\bipartBC}}. \]

These properties extend easily by replacing $\bipartBC$ with $\bipartB$ or $\bipartC$. This gives us an eigenmatrix $P^{\bipartB}$ satisfying
\begin{equation}\label{eq:ebf919}
  X_i = \sum_{j=0}^{t_{\bipartB}} P_{j,i}^{\bipartB} L_j
\end{equation}
and an eigenmatrix $P^{\bipartC}$ satisfying
\[ Y_i = \sum_{j=0}^{t_{\bipartC}} P_{j,i}^{\bipartC} R_j. \]
We also have a dual eigenmatrix $Q^{\bipartB}$ satisfying
\begin{equation}\label{eq:a17266}
L_i = \frac{1}{\abs{\bipartB}} \sum_{j=0}^{t_{\bipartB}} Q_{j,i}^{\bipartB} X_j.
\end{equation}
and $Q^{\bipartC}$ satisfying
\[ R_i = \frac{1}{\abs{\bipartC}} \sum_{j=0}^{t_{\bipartBC}} Q_{j,i}^{\bipartC} Y_j. \]

The (dual) eigenmatrices for \( \bipartB \) and \( \bipartC \) are the (dual) eigenmatrices of symmetric association schemes, and so they have stronger properties that can be found in references such as Chapter~2 of Bannai and Ito~\cite{bannai1984algebraic} or Section~2.~2 of Brouwer, Cohen, and Neumaier~\cite{bcn}.

Of particular note, we can define the multiplicities of the eigenvalues relative to one cell of the partition. Rewriting Equation~\ref{sii826}, we have \( m_r^{\bipartB} = \tr \lp L_r \rp \) and \( m_r^{\bipartC} = \tr \lp R_r \rp. \) Note that
\[ L_i \circ X_0 = \frac{1}{\abs{\bipartB}} Q_{0,i}^{\bipartB} X_0, \]
and taking the sum of every entry on both sides gives us \( m_i^{\bipartB} = Q_{0,i}^{\bipartB}. \) An analogous argument gives us \( m_i^{\bipartC} = Q_{0,i}^{\bipartC}. \)

If \( 0 \leq h \leq \tilde{t}_{\bipartBC}, \) then~\ref{dbc4} tells us that \( L_h = D_h D_h^T \), from which it follows that
\[ m_h^{\bipartB} = \tr \lp L_h \rp = \tr \lp D_h D_h^T \rp = \mathrm{sum} \lp D_h \circ D_h \rp. \]
Using Equation~\ref{eq:68c2ad}, we have
\[ D_h \circ D_h = \frac{1}{\abs{\bipartB} \abs{\bipartC}} \sum_{i=1}^{t_{\bipartBC}} Q_{i,h}^{\bipartBC} Q_{i,h}^{\bipartBC} N_i, \] 
so
\begin{equation}
m_h^{\bipartB} = \mathrm{sum} \lp D_h \circ D_h \rp = \frac{1}{\abs{\bipartC}} \sum_{i=1}^{t_{\bipartBC}} \lp Q_{i,h}^{\bipartBC} \rp^2 k_i^{\bipartBC}.\label{eq:e3efd2}
\end{equation}
Further, since~\ref{dbc4} also tells us that \( R_h = D_h^T D_h, \) we can apply Lemma~\ref{BBT} to see that \( m_h^{\bipartB} = m_h^{\bipartC}. \)
%\begin{equation}\label{eq:mhefd2}
% m_h^{\bipartC} = m_h^{\bipartB} = \frac{1}{\abs{\bipartC}} \sum_{i=0}^{t_{\bipartBC}} \lp Q_{i,h}^{\bipartBC} \rp^2 k_i^{\bipartBC}.
%\end{equation}

\section{Intersection Numbers}\label{secIntersect}

Consider the intersection numbers of a bipartite coherent configuration. To avoid an excess of superscripts, we will use a similar ordering of the matrices as we used in Section~\ref{secPPoly}. This allows us to represent the intersection numbers as two sequences rather than eight, which has particular relevance when considering distance-biregular graphs.

By~\ref{bcc5}, we know there exist coefficients, called intersection numbers, such that for \( 0 \leq i \leq t_{\bipartB} \) and \( 1 \leq j \leq t_{\bipartBC} \) we have
\begin{equation}\label{eq:bbba23}
X_i N_j = \sum_{h=1}^{t_{\bipartBC}} \xi_{2i,2j-1}^{2h-1} N_h.
\end{equation}
Thus for  \( 1 \leq h \leq t_{\bipartBC}, \) we have
\[ N_{h} \circ X_iN_j = \xi_{2i,2j-1}^{2h-1} N_{h}.  \]
Taking the sum of all the entries of both sides, we have
\[ \tr \lp N_{h}^T  X_i N_j \rp= \mathrm{sum} \lp N_{h} \circ X_i N_j \rp = \xi_{2i,2j-1}^{2 h-1} \abs{\bipartB} k_{h}^{\bipartBC}. \]
From Equations~\ref{eq:1aecc5} and~\ref{eq:ebf919}, we have
\[ N_{h}^T X_i N_j = \lp \sum_{r=0}^{\tilde{t}_{\bipartBC}} P_{r,h}^{\bipartBC} D_r^T \rp \lp \sum_{s=0}^{t_{\bipartB}} P_{s,i}^{\bipartB} L_s \rp \lp \sum_{m=0}^{\tilde{t}_{\bipartBC}} P_{m,j}^{\bipartBC} D_m \rp. \]
By~\ref{dbc4}
\[ D_r^T L_s D_m  = D_r^T D_s D_s^T D_m = \delta_{r,s} \delta_{r,m} L_r. \]
Thus
\[ N_{h}^T X_i N_j = \sum_{r=0}^{\tilde{t}_{\bipartBC}} P_{r, h}^{\bipartBC} P_{r,i}^{\bipartB} P_{r,j}^{\bipartBC}L_r, \]
so
\[ \xi_{2i,2j-1}^{2 h-1} = \frac{1}{\abs{\bipartB} k_{h}^{\bipartBC}} \sum_{r=0}^{\tilde{t}_{\bipartBC}} m_r^{\bipartB} P_{r, h}^{\bipartBC} P_{r,i}^{\bipartB} P_{r,j}^{\bipartBC}. \]
It follows that $\xi_{2i,2j+1}^{2 h+1}$ is determined by the eigenmatrices $P^{\bipartB}$ and \( P^{\bipartBC}. \)

Similarly,
\[ \xi_{2i-1,2j}^{2 h-1} = \frac{1}{\mathrm{sum} \lp N_h \rp} \mathrm{sum} \lp N_h \circ N_i Y_j \rp = \frac{1}{\abs{\bipartB} k_{h}^{\bipartBC}} \sum_{r=0}^{\tilde{t}_{\bipartBC}} m_r^{\bipartC} P_{r,h}^{\bipartBC} P_{r,i}^{\bipartBC} P_{r,j}^{\bipartC}; \]
\begin{equation}\label{eq:ae1d52} \xi_{2i,2j}^{2h} = \frac{1}{\mathrm{sum} \lp X_h \rp} \mathrm{sum} \lp X_{h} \circ X_i X_j^T \rp = \frac{1}{\abs{\bipartB} k_{h}^{\bipartB}}\sum_{r=0}^{t_{\bipartB}} m_r^{\bipartB} P_{r,h}^{\bipartB} P_{r,i}^{\bipartB} P_{r,j}^{\bipartB};
\end{equation}
and
\[ \xi_{2i-1,2j-1}^{2h} = \frac{1}{\mathrm{sum} \lp X_{h} \rp} \mathrm{sum} \lp X_{h} \circ N_i N_j^T \rp = \frac{1}{\abs{\bipartB} k_{h}^{\bipartB}} \sum_{r=0}^{\tilde{t}_{\bipartBC}} m_r^{\bipartB} P_{r,h+1}^{\bipartB} P_{r,i}^{\bipartBC} P_{r,j}^{\bipartBC}. \]
Note \( \xi_{i,j}^{h} = 0 \) for all \( h > t_{\bipartB} \) if \( h \) is even and \( h > t_{\bipartBC} \) if \( h \) is odd.

Flipping \( \bipartB, \bipartC, \) and the relevant eigenmatrices gives us the intersection numbers for \( \bipartC, \) which we denote by \( \sigma_{i,j}^{h}. \)

We can think of the \( X_i \) as the adjacency matrices of graphs on vertex set \( \bipartB, \) and \( N_i \) as the biadjacency matrix from \( \bipartB \) to \( \bipartC \). With some abuse of notation, we denote the graphs by the matrices.

Fix a vertex \( \vtxa \in \bipartB. \) For \( 1 \leq h \leq t_{\bipartBC}, \) we choose a vertex \( \vtxb \) such that \( \vtxa \) and \( \vtxb \) are adjacent in \( N_h \). Then for \( 0 \leq j \leq t_{\bipartBC}, \) we use Equation~\ref{eq:bbba23} to see that
\[ \lp X_i N_j \rp_{\vtxa, \vtxb} = \xi_{2i,2j-1}^{2h-1}. \]
We also have
\[ \lp X_i N_j \rp_{\vtxa, \vtxb} = \sum_{\vtxc \in \bipartB} \lp X_i \rp_{\vtxa, \vtxc} \lp N_j \rp_{\vtxc, \vtxb} = \abs{\lsb \vtxc \in \bipartB: \vtxa \sim_{X_i} \vtxc, \vtxb \sim_{N_j} \vtxc \rsb}, \]
so \( \xi_{2i,2j-1}^{2h-1} \) must count some discrete number of vertices, and as such, must be a non-negative integer. Similar arguments tell us that all \( \xi_{i,j}^{h}, \sigma_{i,j}^h \) must be non-negative integers.

The intersection numbers of distance-biregular graphs was considered in~\cite{PhD}, especially in Chapters~3 and~6.

\section{Krein Parameters}\label{secKrein}

The dual of intersection numbers are the Krein parameters. Krein parameters were defined for distance-biregular graphs by Delorme~\cite{delorme} and for coherent configurations by Hobart~\cite{hobart1995krein} and Ito and Munemasa~\cite{ito2020krein}. We denote \( \rho_{rs} \lp t \rp, \lambda_{rs} \lp t \rp \) and \( \Delta_{r,s} \lp t \rp \) as the coefficients such that
\[ L_i \circ L_j = \frac{1}{\abs{\bipartB}} \sum_{h=0}^{t_{\bipartB}} \lambda_{ij} \lp h \rp L_h; \]
\[ D_i \circ D_j = \frac{1}{\sqrt{\abs{\bipartB} \abs{\bipartC}}} \sum_{h=0}^{t_{\bipartBC}} \Delta_{ij} \lp h \rp D_h; \]
and
\[ R_i \circ R_j = \frac{1}{\abs{\bipartC}} \sum_{h=0}^{t_{\bipartC}} \rho_{ij} \lp h \rp R_h. \]

For \( 0 \leq i, j, h \leq \tilde{t}_{\bipartBC}, \) we have that
\[ D_h^T \lp D_i \circ D_j \rp = \frac{1}{\sqrt{\abs{\bipartB} \abs{\bipartC}}} \Delta_{ij} \lp h \rp L_h. \]
Summing over the $\lp \vtxa, \vtxa \rp$ entry of both sides for all $\vtxa \in \bipartB$, we see
\[ \Delta_{ij} \lp h \rp = \frac{\sqrt{\abs{\bipartB} \abs{\bipartC}}}{m_h^{\bipartB}} \mathrm{tr} \lp D_h^T \lp D_i \circ D_j \rp \rp = \frac{\sqrt{\abs{\bipartB} \abs{\bipartC}}}{m_h^{\bipartB}} \mathrm{sum} \lp D_h \circ D_i \circ D_j \rp. \]
Using Equation~\ref{eq:68c2ad}, we have
\[ D_h \circ D_i \circ D_j = \frac{1}{\sqrt{\abs{\bipartB} \abs{\bipartC}}^3} \sum_{\ell=0}^{\tilde{t}_{\bipartBC}} Q^{\bipartBC}_{\ell,h} Q^{\bipartBC}_{\ell,i} Q^{\bipartBC}_{\ell,j} N_{\ell}, \]
so
\begin{equation}
\Delta_{ij} \lp h \rp = \frac{1}{\abs{\bipartC} m_h^{\bipartB}} \sum_{\ell=0}^{\tilde{t}_{\bipartBC}} Q^{\bipartBC}_{\ell,h} Q^{\bipartBC}_{\ell,i} Q^{\bipartBC}_{\ell,j} k_{\ell}^{\bipartBC}.\label{eq:1880dd}
\end{equation}
Recall that \( m_h^{\bipartB} = m_h^{\bipartC}, \) and from Equation~\ref{eq:kl2014} we have that
\[ \frac{k_{\ell}^{\bipartBC}}{\abs{\bipartC}} = \frac{1}{\sqrt{\abs{\bipartB} \abs{\bipartC}}} P_{0,\ell}^{\bipartBC} = \frac{k_{\ell}^{\bipartCB}}{\abs{\bipartB}} \]
so Equation~\ref{eq:1880dd} does not depend on the choice of cell to sum over.

Similar arguments give us
\begin{equation}
\lambda_{ij} \lp h \rp = \frac{1}{\abs{\bipartB} m_h^{\bipartB}} \sum_{\ell=0}^{t_{\bipartB}} Q_{\ell,h}^{\bipartB} Q_{\ell,i}^{\bipartB} Q_{\ell,j}^{\bipartB} k^{\bipartB}_{\ell}\label{eq:9301ea}
\end{equation}
and
\begin{equation}
\rho_{ij} \lp h \rp = \frac{1}{\abs{\bipartC} m_h^{\bipartC}} \sum_{\ell=0}^{t_{\bipartC}} Q_{\ell,h}^{\bipartC} Q_{\ell,i}^{\bipartC} Q_{\ell,j}^{\bipartC} k^{\bipartC}_{\ell}.\label{eq:94b263}
\end{equation}
Equations~\ref{eq:9301ea} and~\ref{eq:94b263} can also be  derived using results from association schemes.

As with the intersection numbers, we can use the Krein parameters to get a necessary condition for a bipartite coherent configuration.

\begin{prop}[Delorme~\cite{delorme}, Hobart~\cite{hobart1995krein}, Ito and Munemasa~\cite{ito2020krein}]\label{kreinInequality}For any bipartite coherent configuration, the Krein parameters satisfy the following:
  \begin{enumerate}[label=(\roman*)]
  \item\label{krein1} For any $0 \leq i, j, h \leq t_{\bipartB}$,
\[ \lambda_{ij} \lp h \rp \geq 0; \]
\item\label{krein2} For any $0 \leq i, j, h \leq t_{\bipartC}$,
\[ \rho_{ij} \lp h \rp \geq 0; \]
\item\label{krein3} For any $0 \leq i, j, h \leq \tilde{t}_{\bipartBC}$,
  \[ \rho_{ij} \lp h \rp \lambda_{ij} \lp h \rp \geq \Delta_{ij} \lp h \rp^2. \]
  \end{enumerate}
\end{prop}

\proof Delorme's proof of the Krein inequality for distance-biregular graphs~\cite{delorme} extends directly to a proof of Proposition~\ref{kreinInequality}. Alternatively, Proposition~\ref{kreinInequality} can be seen as a special case of Krein conditions on a larger classes of coherent configurations. Hobart~\cite{hobart1995krein} phrased Krein conditions on coherent configurations in terms of certain matrices being positive semidefinite, and Ito and Munemasa~\cite{ito2020krein} proved that for fibre-commutative coherent configurations, it was sufficient to prove the Krein matrix
\[ \pmat{\lambda_{ij} \lp h \rp & \Delta_{ij} \lp h \rp\\ \Delta_{ij} \lp h \rp & \rho_{ij} \lp h \rp} \]
is positive semidefinite. Equivalently, the principal minors are all non-negative, which gives us the desired inequalities.\qed

Delorme~\cite{delorme} and Hobart~\cite{hobart1995krein} both provided examples of bipartite coherent configurations where condition~\ref{krein3} failed, even though conditions~\ref{krein1} and~\ref{krein2} held. By contrast, Ito and Munemasa~\cite{ito2020krein} proved that for generalized quadrangles, the Krein parameters of the bipartite coherent configuration reduced to~\ref{krein1} and~\ref{krein2}.

The inequalities in Proposition~\ref{kreinInequality} had been worked out for some special cases of bipartite coherent configurations prior to the work of Delorme~\cite{delorme}, Hobart~\cite{hobart1995krein} and Ito and Munemasa~\cite{ito2020krein}. Krein conditions were used to restrict parameters of generalized polygons~\cite{haemers1981inequality, higman1975invariant}, and the inequalities obtained for quasi-symmetric designs~\cite{calderbank1988inequalities,neumaier1982regular} matched those found by Hobart~\cite{hobart1995krein} when considering the Krein condition for coherent configurations.

The work of Hobart~\cite{hobart1995krein} and Ito and Munemasa~\cite{ito2020krein} applies in the more general case of coherent configurations with an arbitrary number of commuting fibres. However, there are no known examples where the Krein conditions succeed on the \( 1 \times 1 \) and \( 2 \times 2 \) principal minors, but fail on larger principal minors. It is also worth noting that the examples studied using Krein parameters of coherent configurations are all bipartite coherent configurations.

\section{$Q$-Polynomial Coherent Configurations and Duality}\label{duality}

Suda~\cite{suda2022q} gave many examples of $Q$-polynomial coherent configurations. An obvious necessary condition for a coherent configuration to be a bipartite coherent configuration is for it to have two fibres, though Example~\ref{hobart} shows that this condition is not sufficient. However, combined with the stronger conditions for the coherent configurations that Suda studied, we can get a sufficient condition for certain types of two-fibre coherent configurations.

\begin{lem}\label{t1tt1}Let \( \mathcal{C} \) be a fibre-symmetric coherent configuration satisfying~\ref{suda1}-~\ref{suda4}. If \( \mathcal{C} \) has type
  \[ \pmat{t+1 & t \\ & t+1} \]
  then it is a bipartite coherent configuration.
\end{lem}

\proof We are dealing with a two-fibre coherent configuration with a distinguished second basis, so we will use the same notation introduced in Sections~\ref{secBCC} and~\ref{secBasis}. Since \( \mathcal{C} \) is a coherent configuration, we know that \( \mathcal{C} \) satisfies~\ref{bcc1}-~\ref{bcc5}. Condition~\ref{bcc4} follows from \( \mathcal{C} \) being fibre-symmetric. Thus it only remains to show~\ref{bcc6}.

Let \( 0 \leq i \leq t \). Conditions~\ref{suda1}-~\ref{suda4} guarantee a second basis, so we write
\begin{equation}\label{ltly4157}
X_i = \sum_{r=0}^{t} P^{\bipartB}_{r,i} L_r.
\end{equation}
Now by~\ref{suda3} and~\ref{suda4}, we know that for all \( 0 \leq r \leq t-1, \) we have \( D_r D_r^T  = L_r, \) so Equation~\ref{ltly4157} becomes
\[ X_i = \sum_{r=0}^{t-1} P^{\bipartB}_{r,i} D_r D_r^T + P^{\bipartB}_{t,i} \lp I - \sum_{r=0}^{t-1} D_r D_r^T \rp = \sum_{r=0}^{t-1} \lp P^{\bipartB}_{r,i} - P^{\bipartB}_{t,i} \rp D_r D_r^T + P^{\bipartB}_{t,i} I. \]
Rewriting \( D_r \) and \( D_r^T \) in the basis of Schur idempotents, we have that
\[ X_i \in \Span \lsb N_j N_h^T : 0 \leq j, h \leq t \rsb \cup \lsb I \rsb. \]
An analogous argument holds for \( Y_i, \) establishing~\ref{bcc5}.\qed

We use this result to sketch one particular example mentioned in Suda~\cite{suda2022q}.

\begin{ex}Beginning with a \( Q \)-polynomial association scheme, Delsarte~\cite{delsarte1977pairs} defined relative $t$-designs. Bannai and Bannai~\cite{bannai2012remarks} proved a Fisher-type inequality, and considered the case of when the bound is tight. A relative $t$-design can be supported on some number of subconstituents, or shells.

  The \( n \)-dimensional hypercube is a $Q$-polynomial association scheme. The tight relative $t$-designs in the hypercube that are supported on one shell are precisely \( t-\lp n, k, \lambda \rp \)-combinatorial designs. Bannai, Bannai, Tanaka, and Zhu~\cite{bannaiBannai} proved that tight relative $2e$-designs supported on two shells gives rise to a coherent configuration of type
  \[ \pmat{e+1 & e \\ & e+1}. \]
  Suda~\cite{suda2022q} showed that this coherent configuration is $Q$-polynomial, and thus an example of a $Q$-polynomial bipartite coherent configuration.
\end{ex}

There is a rich theory connecting design theory to $Q$-polynomial association schemes, beginning with the work of Delsarte~\cite{delsarte1973algebraic}. Many of the examples of Suda~\cite{suda2022q} build on these connections by constructing $Q$-polynomial coherent configurations from the designs related to $Q$-polynomial association schemes, and one of the open problems at the end of that paper was to develop design theory in the context of $Q$-polynomial coherent configurations. A special case would be to develop design theory in the context of $Q$-polynomial bipartite coherent configurations.

There is another, more direct, connection between design theory and bipartite coherent configurations. A bipartite coherent configuration can be seen as modelling points, blocks, and various incidence relations with particularly strong structure, such as quasi-symmetric designs of Example~\ref{quasisymmetric} or the strongly regular designs of Example~\ref{stronglyRegularDesign}. This makes bipartite coherent configurations a logical place to extend some of the theory of association schemes.

\begin{prob}Can we develop design theory in the context of bipartite coherent configurations?
\end{prob}

The connections between design theory and $Q$-polynomial association schemes is only one part of the interest of $Q$-polynomial association schemes. There is a well-developed theory of duality between $P$- and $Q$-polynomial association schemes, which has been discussed in Bannai and Ito~\cite{bannai1984algebraic}, Brouwer, Cohen, Neumaier~\cite{bcn}, and the more recent survey by Martin and Tanaka~\cite{martin2009commutative}. Bipartite coherent configurations give a framework to develop similar duality for distance-biregular graphs.

\begin{prob}In the context of bipartite coherent configurations, what are other equivalent definitions of \( Q \)-polynomial?
\end{prob}

In particular, for association schemes we can define a notion of \textit{formally dual} (see, for instance, Chapter~12 of Godsil~\cite{blue}) and any association scheme that is formally dual to a \( P \)-polynomial association scheme is \( Q \)-polynomial.

There is a beautiful theory of $Q$-polynomial distance-regular graphs. Of particular note, Leonard~\cite{leonard1982orthogonal,leonard1984metric,leonard1984parameters} characterized \( Q \)-polynomial distance-regular graphs and showed that they could be expressed in terms of a small number of parameters. A number of more recent results on \( Q \)-polynomial distance-regular graphs can be found in the survey by Terwilliger~\cite{terwilliger2022distance}.

\begin{prob}Extend the theory of \( Q \)-polynomial distance-regular graphs to \( Q \)-polynomial distance-biregular graphs.
\end{prob}

Finally, it would be interesting to have more examples of bipartite coherent configurations, in and outside the context of \( P \)- and \( Q \)-polynomial bipartite coherent configurations.

\section*{Acknowledgements}
I would like to thank Soff\'ia \'Arnad\'ottir for feedback on previous drafts, and the anonymous referees whose feedback considerably improved the presentation. This work was supported by Kempe foundation JCSMK22-0160.

%\bibliographystyle{acm}
%\bibliography{/home/sabrina/Documents/Academics/Waterloo/PhD/ExtremeGhostbusting/Bibliography}

\end{document}